\documentclass[conference]{IEEEtran}
\IEEEoverridecommandlockouts

\usepackage{cite}
\usepackage{amsmath,amssymb,amsfonts}
\usepackage{acronym}
\usepackage{graphicx}
\usepackage{textcomp}
\usepackage{xcolor}
\usepackage{threeparttable}
\usepackage{booktabs}

\usepackage{pifont}
\usepackage{caption}
\usepackage{subcaption}
\usepackage[hidelinks]{hyperref}
\usepackage{comment}
\usepackage{upgreek}

\usepackage{multirow}

\def\BibTeX{{\rm B\kern-.05em{\sc i\kern-.025em b}\kern-.08em
T\kern-.1667em\lower.7ex\hbox{E}\kern-.125emX}}

\newcommand{\bfb}[1]{\boldsymbol{\rm #1}}
\newcommand{\bfg}[1]{\boldsymbol{#1}}
\newcommand{\nx}{\nu}
\newcommand{\ny}{\mu}

\newcommand{\ntdi}{q}
\newcommand{\xpc}{\xi}

\newcommand{\reigv}{\bfg v}
\newcommand{\reigvel}{v}
\newcommand{\reigvmat}{\bfg U}
\newcommand{\leigv}{\bfg w}
\newcommand{\leigvel}{w}
\newcommand{\leigvmat}{\bfg W}
\newcommand{\jac}[2]{\bfg{{#1}}_{\hspace{-0.2mm}#2}}
\newcommand{\jj}{\jmath}
\newcommand{\zt}{\hat z}
\newcommand{\T}{^{\intercal}}
\newcommand{\AS}{\bfb A}

\newcommand{\Atdi}{\bfb{G}}
\newcommand{\reigvtdi}{\bfb{v}}
\newcommand{\leigvtdi}{\bfb{w}}
\newcommand{\leigvmattdi}{\boldsymbol{\mathcal{W}}}
\newcommand{\reigvmattdi}{\boldsymbol{\mathcal{V}}}

\newcommand{\PF}{{p}}
\newcommand{\PFmat}{{\bfb P}}
\newcommand{\PFmattdi}{\bfb{\Pi}}

\newcommand{\ci}{c}
\newcommand{\cii}{b}

\newcommand{\errs}{\epsilon_s}
\newcommand{\errp}{\epsilon_p}

\newcommand{\dirka}{\alpha}
\newcommand{\dirkb}{\beta}

\acrodef{fem}[FEM]{forward Euler method}
\acrodef{hm}[HM]{Heun's method}
\acrodef{tm}[TM]{trapezoidal method}
\acrodef{pf}[PF]{participation factor}
\acrodef{sssa}[SSSA]{small-signal stability analysis}
\acrodef{bem}[BEM]{backward Euler method}
\acrodef{dirk}[2S-DIRK]{two-stage diagonally implicit Runge Kutta}
\acrodef{dae}[DAE]{differential-algebraic equation}
\acrodef{pll}[PLL]{phase-locked loop}
\acrodef{sg}[SG]{synchronous generator}
\acrodef{der}[DER]{distributed energy resource}

\acrodef{ode}[ODE]{ordinary differential equation}
\acrodef{tdi}[TDI]{time-domain integration}
\acrodef{agc}[AGC]{automatic generation control}
\acrodef{tg}[TG]{turbine governor}

\begin{document}

\title{Mode-Shape Deformation of Power System
DAEs by Time-Domain Integration Methods
\thanks{G. Tzounas and G.~Hug are supported by the Swiss National Science Foundation under project NCCR Automation (grant no.~51NF40 18054); C.~Tajoli and G.~Hug are supported by project ReMaP.}
} 

\author{\IEEEauthorblockN{Carlo Tajoli, {\em IEEE Student Member}, Georgios Tzounas, {\em IEEE Member}, and Gabriela Hug, {\em IEEE Senior Member}}
\IEEEauthorblockA{\textit{Power Systems Laboratory} \\
\textit{ETH Z{\"u}rich}\\
Zurich, Switzerland \\
\{tajoli, georgios.tzounas, hug\}@eeh.ee.ethz.ch}
}

\maketitle 

\begin{abstract} 
This paper studies the numerical deformation that 
\ac{tdi} methods introduce to the 
shape of the coupling between the dynamic modes and variables of power system models.
To this aim, we employ a \ac{sssa}-based framework where such \textit{mode-shape deformation} is efficiently identified by comparing the modal \acp{pf} of the power system model with the \acp{pf} of the discrete-time system that is derived from the application of the \ac{tdi} method.  The proposed approach is illustrated for \ac{tdi} methods commonly used in dynamic power system calculations.  
\end{abstract}


\section{Introduction}
\label{sec:intro}

\subsection{Motivation}

The stability analysis of a power system following a large disturbance -- such as the sudden loss of an important generator,
a 
line fault, etc. -- relies on the 
solution of a non-linear model of \acp{dae} 
\cite{kundur:94}.  
Power system software tools approximate this solution numerically
by running a time-domain simulation routine.  However, rapid and precise stability analysis through time-domain simulations is 
not straightforward, especially with the growing penetration of converter-based resources which significantly increases the dynamic complexity and stiffness of power system models.

\subsection{Literature Review}

There exist two 
time-domain simulation approaches to obtain the solution of the \acp{dae} that describe the dynamics of power systems, namely simultaneous and partitioned \cite{machowski2020power}.  
In the simultaneous approach, differential and algebraic equations are solved together as one set at each time step through an implicit integration method, such as the Theta method \cite{1995paserba,powerfactory}.
In the partitioned approach, on the other hand, differential equations are solved at each step for state variables, whereas algebraic equations are solved separately. The solution of differential equations in this case is typically obtained with an explicit integration method \cite{stott:1979}.  For example, a family of methods commonly employed in a partitioned-solution setup is that of explicit Adams-Bashforth \cite{psse_pag2,ge:pslf}.

In contrast to implicit \acf{tdi} methods, explicit methods are known to be prone to numerical instabilities.  This limits the ability of these methods to use large integration time steps and has often driven efforts for the development of device models that are numerically robust when combined with a given commercial explicit solver.  In this vein, recent works have focused on the formulation and numerical robustness of converter-based resource models for systems with low short-circuit strength, e.g.,~see \cite{ramasubramanian2020positive,ramasubramanian2022parameterization}.

The accuracy of a time-domain simulation is traditionally evaluated through truncation error analysis.  Moreover, the numerical stability of a \ac{tdi} method is conventionally characterized by testing its convergence on a linear scalar equation.  Recent studies focused on the development of a framework to assess accuracy and numerical stability of \ac{tdi} methods in a unified way.  In particular, \cite{tzounas2022tdistab,tzounas2022tdistabd,tzounas:psastab} estimate the numerical distortion that a given \ac{tdi} method introduces to a power system model by comparing the small-signal dynamic modes of the original model
with the modes of the approximated system that results from the application of the method.  Such framework allows, first, to extract useful upper time step bounds that satisfy prescribed requirements of precision and over/under-damping; and, second, to provide a fair computational comparison among different methods. 

Apart from the numerical error that they cause to the dynamic modes of a power system model, \ac{tdi} methods may also introduce a spurious deformation to the shape of the coupling between dynamic modes and system variables.  This \textit{mode-shape} aspect of numerical deformation
has not, to the best of our knowledge, been investigated in the literature. To provide a first study that 
tackles this aspect is the main goal of our work in this paper.

\subsection{Contributions}

The contributions of the paper are twofold, as follows:

\begin{itemize}
    \item  Provision of a \acf{sssa}-based technique to estimate the numerical mode-shape deformation introduced to \ac{dae} power system models by \ac{tdi} methods.
    \item Thorough discussion on the mode-shape deformation caused to power system \acp{dae} by well-known \ac{tdi} methods, including 
    Theta, \ac{dirk}, 
    and \ac{hm}.
    
\end{itemize}

\subsection{Paper Organization}

The remainder of the paper is organized as follows.  Section~\ref{sec:tdi} recalls the formulation and numerical integration of \ac{dae} power system models.  Section~\ref{sec:numdeform} describes the proposed approach to quantify the numerical deformation of the system's mode shapes caused by \ac{tdi} methods. Section~\ref{sec:case} discusses the case study.  Conclusions are drawn in Section~\ref{sec:conclusion}.

\section{Power System Model and Numerical Solution}
\label{sec:tdi}

\subsection{DAE Model}
\label{sec:tdi:model}

In short-term stability analysis, the dynamic model of a power system is conventionally formulated as a set of non-linear \acp{dae}, as follows \cite{kundur:94}:
\begin{equation}
  \begin{aligned}
    \label{eq:dae}
    {\bfg x}'(t) &= \bfg f(\bfg x(t), \bfg y (t) 
    ) \, , \\
    \bfg 0_{\ny,1} &=  \bfg g(\bfg x(t), \bfg y(t)
    ) \, .
  \end{aligned}    
\end{equation}
In \eqref{eq:dae}, $\bfg x(t) : [0,\infty) \rightarrow \mathbb{R}^{\nx}$ and
$\bfg y(t): [0,\infty) \rightarrow \mathbb{R}^{\ny}$ are the states and algebraic variables, respectively, of the system; 
$\bfg f : \mathbb{R}^{\nx+\ny} \rightarrow \mathbb{R}^{\nx}$ and $\bfg g : \mathbb{R}^{\nx+\ny} \rightarrow \mathbb{R}^{\ny}$ are non-linear functions; $\bfg 0_{\ny,1}$ is the zero matrix of dimensions $\ny \times 1$.  For simplicity, discrete system dynamics are not explicitly considered in \eqref{eq:dae}. Readers interested in the modeling and handling of discontinuities are referred to \cite{2022ampsas} and the bibliography therein.

\subsection{Numerical Integration}

The time-domain simulation of a power system model consists in employing a proper numerical method to approximate the solution of \eqref{eq:dae} for a known set of initial conditions.  Every numerical \ac{tdi} method applied to \eqref{eq:dae} can be mathematically described as a set of non-linear difference equations whose definition depends on $\bfg f$ and $\bfg g$.  For example, employing the well-known
Theta method \cite{1995paserba} leads to the following set of difference equations:
\begin{equation}
\begin{aligned}\label{eq:theta}
  \bfg x_{n+1} &= \bfg x_{n} + h
 [\theta\bfg f(\bfg x_{n},\bfg y_{n}) + 
 (1-\theta)\bfg f(\bfg x_{n+1},\bfg y_{n+1})] \, , \\
 \bfg 0_{\ny,1} &= h \bfg g(\bfg x_{n+1},\bfg y_{n+1}) \, ,
\end{aligned}
\end{equation}
where $h$ is the simulation time step size; $0\leq\theta\leq0.5$ 
determines the method's damping;
and $\bfg x_{n+1-\ell} : \mathbb{N} \rightarrow\mathbb{R}^{\nx}$, $\bfg y_{n+1-\ell} : \mathbb{N} \rightarrow \mathbb{R}^{\ny}$,
$l =\{0,1\}$.
Given the values of state and algebraic variables at some point $(\bfg x_{n}, \bfg y_{n}) := [\bfg x_{n}\T, \bfg y_{n}\T]\T$ (where $\T$ is the matrix transpose), the goal at each time step is to compute the new values $(\bfg x_{n+1}, \bfg y_{n+1})$.  The latter provides an approximation of the exact solution of \eqref{eq:dae},~i.e:
\begin{equation}
\begin{aligned}
\label{eq:discrxy}
\bfg x_{n+1-\ell} &\approx \bfg x(t+(1-\ell) h) \, , \\
\bfg y_{n+1-\ell} &\approx \bfg y(t+(1-\ell) h) \, .
\end{aligned}    
\end{equation}
The accuracy and convergence of a \ac{tdi} 
depends on the time step size $h$, as well as on the numerical properties of the \ac{tdi} method employed. For example, \eqref{eq:theta} for $\theta=0.5$ corresponds to the \ac{tm}, which 
always converges for stable and diverges for unstable trajectories.
On the other hand, for $\theta=0$, \eqref{eq:theta} corresponds to the \ac{bem}, which has very fast convergence but tends to overdamp the dynamics of the system.

\subsection{Model Stiffness and SSSA}
\label{sec:sssa}

System \eqref{eq:dae} is known to be stiff, i.e.~its differential equations span a wide range of time constants \cite{kundur:94}.
The stiffness of \eqref{eq:dae} can be measured by the ratio between the largest and smallest
eigenvalues of the corresponding small-signal model.

Consider that a stationary solution $(\bfg x_o,\bfg y_o)$ of \eqref{eq:dae} is known.
Then, differentiating \eqref{eq:dae} at the stationary point gives:
\begin{equation}
  \begin{aligned}
    \label{eq:dae:lin}
    \tilde {\bfg x}'(t)  &= 
    \jac{f}{x} \tilde {\bfg x}(t) + \jac{f}{y} \tilde {\bfg y}(t) \, , \\
    \bfg 0_{\ny,1} &=\jac{g}{x} \tilde {\bfg x}(t) + \jac{g}{y} \tilde {\bfg y}(t) \, , 
  \end{aligned}    
\end{equation}
where $\tilde {\bfg x}(t) = \bfg x(t) - \bfg x_o$, $\tilde {\bfg y}(t) = \bfg y(t) - \bfg y_o$; and $\jac{f}{x}$, $\jac{f}{y}$, $\jac{g}{x}$, $\jac{g}{y}$ are Jacobian matrices evaluated at $(\bfg x_o, \bfg y_o)$.  Under the assumption that $\bfg g_y$ is non-singular\footnote{In this paper, we assume that $\bfg g_y$ is invertible. This 
assumption 
comes with no loss of generality, as
potential singularities of $\bfg g_y$ can be always eliminated by reformulating \eqref{eq:dae} to an equivalent \ac{dae} set with non-singular $\bfg g_y$.}, algebraic variables can be eliminated\footnote{Eliminating $\tilde {\bfg y}$ is the best approach for small/medium size systems. In large systems 
it is more efficient to maintain sparsity 
and work directly with \eqref{eq:dae:lin}. } and \eqref{eq:dae:lin} can be rewritten as a set of linear ordinary differential equations, as follows:
\begin{equation}
  \begin{aligned}
    \label{eq:lin:ode}
    \tilde {\bfg x}'(t)  &= \AS  
     \tilde {\bfg x}(t)  \, , 
  \end{aligned}    
\end{equation}
where $\AS =\jac{f}{x}-\jac{f}{y}\jac{g}{y}^{-1}\jac{g}{x}$. 
Then, stability of \eqref{eq:dae:lin} is assessed through the eigenvalues of \eqref{eq:lin:ode}, which are obtained from the numerical solution of the algebraic problem \cite{book:eigenvalue}:
\begin{align}\label{eq:ep}
(s\bfb I_\nx -\AS) \reigv &= \bfg 0_{\nx,1}  \, , \\
\leigv (s\bfb I_\nx-\AS) &= \bfg 0_{1,\nx}   \, ,
\label{eq:rep}
\end{align}
where $s$ denotes a complex frequency in the $S$-domain; $\bfb I_\nx$ denotes the identity matrix of dimensions $\nx \times \nx$;  $\reigv \in \mathbb{C}^{\nx\times 1}$ and $\leigv \in \mathbb{C}^{1\times \nx}$.  
Every $s_i$, $i=\{1,2,\ldots,\nx\}$, that satisfies \eqref{eq:ep} is an eigenvalue of $\AS$,
with $\reigv_i$, $\leigv_i$ being the corresponding right and left, respectively, eigenvectors.  Then, the system is asymptotically stable if $\forall s_i$, $\Re(s_i)<0$.
Let the system be stable and $s^{\max}$, $s^{\min}$ be the eigenvalues with largest and smallest magnitudes, i.e.~$s^{\max}=\max{|s_i|}$,
$s^{\min}=\min{|s_i|}$, $\forall s_i$, then the stiffness ratio of \eqref{eq:dae} can be defined as follows:
\begin{equation}\label{eq:stiffness}
\mathcal{S} = |s^{\max}|/|s^{\min}| \, .
\end{equation}

\section{Proposed Approach}
\label{sec:numdeform}

\subsection{SSSA of Integration Methods}

The small-disturbance properties of a 
\ac{tdi} method applied to a power system model
can be seen by studying a linear system of difference equations in the form  \cite{tzounas2022tdistab,tzounas:psastab}:
\begin{equation}\label{eq:difference}
\bfb y_{n+1} = \Atdi \bfb y_{n} \, ,
\end{equation}
where $\bfb y_{n}: \mathbb{N} \rightarrow \mathbb{R}^{\ntdi}$.  Equation
\eqref{eq:difference} is a discrete-time approximation of \eqref{eq:lin:ode}, where 
$\Atdi$ varies for different \ac{tdi} methods but is always a function of $\AS$ and $h$.
For the sake of example, consider the Theta method described by \eqref{eq:theta}.  Differentiating \eqref{eq:theta} at $(\bfg x_o,\bfg y_o)$ gives:
\begin{align}
\nonumber
\tilde {\bfg x}_{n+1} = \tilde {\bfg x}_{n} &+ h
 [\theta(\jac{f}{x}\tilde{\bfg x}_{n}+ \jac{f}{y}\tilde{\bfg y}_{n})   \\ 
 \label{eq:theta:lin}
 &+
 (1-\theta)(\jac{f}{x}\tilde{\bfg x}_{n+1} +
  \jac{f}{y}\tilde{\bfg y}_{n+1})] \, , \\
  \label{eq:theta:lin:gn}
\bfg 0_{\ny,1} = \jac{g}{x} & \tilde {\bfg x}_{n+1}
+
 \jac{g}{y} \tilde {\bfg y}_{n+1} \, .
\end{align}
From \eqref{eq:theta:lin:gn} we have that $\tilde{\bfg y}_{n+1} = -\jac{g}{y}^{-1}\jac{g}{x}\tilde {\bfg x}_{n+1}$ and $\tilde{\bfg y}_{n} = -\jac{g}{y}^{-1}\jac{g}{x}\tilde{\bfg x}_{n}$, and 
\eqref{eq:theta:lin}-\eqref{eq:theta:lin:gn}
can be rewritten as follows:
\begin{align}
\label{eq:theta:lin2}
\left[\bfb I_\nx - h(1-\theta)\AS
\right ]\tilde {\bfg x}_{n+1} &= 
(\bfb I_\nx + h\theta\AS) 
\tilde {\bfg x}_{n}  \, ,
\end{align}
or equivalently, 
\begin{align}
\label{eq:theta:lin3}
\tilde {\bfg x}_{n+1} &= 
\left[\bfb I_\nx - h(1-\theta)\AS
\right ]^{-1}
(\bfb I_\nx + h\theta\AS) 
\tilde {\bfg x}_{n}  \, ,
\end{align}
which is a system in the form of \eqref{eq:difference}, where $\bfb y_{n} \equiv \tilde {\bfg x}_{n}$, and:
\begin{equation}
  \begin{aligned}
    \label{eq:theta:pencil}
    \Atdi = \left[\bfb I_\nx - h(1-\theta)\AS
\right ]^{-1}(\bfb I_\nx + h\theta\AS) \, .
  \end{aligned} 
\end{equation}
%
%
The eigenvalue problem associated to \eqref{eq:difference} is:
\begin{align}\label{eq:ep:tdi}
(\zt \bfb I_\nx-\Atdi) \reigvtdi &= \bfg 0_{\ntdi,1} \, , \\
\leigvtdi (\zt\bfb I_\nx -\Atdi) &= \bfg 0_{1,\ntdi} \, ,
\label{eq:rep:tdi}
\end{align}
where $\zt$ is a complex frequency in the $Z$-domain. Then, \eqref{eq:difference} is asymptotically stable if and only if 
$|\zt_j|<1$ $\forall \zt_j$, $j={1,2,\ldots,q}$ 
that satisfies \eqref{eq:ep:tdi},
\eqref{eq:rep:tdi}.
%
Comparison of the eigenvalues of 
$\Atdi$ and $\AS$ provides a rough yet accurate estimate of the numerical deformation that a given \ac{tdi} method introduces when applied to \eqref{eq:dae} \cite{tzounas2022tdistab,tzounas:psastab}.  Obviously, for the eigenvalues of the two matrices to be comparable, they need to be referred to the same plane through the map $z = e^{sh}$.  Let $s_i$ be an eigenvalue of $\AS$ and $\zt_j$ be the corresponding eigenvalue as deformed by the \ac{tdi} method. Then, the associated numerical deformation can be estimated through the relative error:
\begin{equation}\label{eq:s:relerr}
\errs = 100 {(|s_i-{\rm log}(\zt_j)/{h}|)}/{|s_i|}  \, .
\end{equation}
%

%

%

\subsection{Deformation of Mode Shapes}

Apart from the numerical deformation that they introduce to the dynamic modes of a model, \ac{tdi} methods may also deform the coupling shape of dynamic modes and state variables.  In this section, we describe the proposed approach to estimate such 
mode-shape deformation. 

In the context of \ac{sssa}, the information of mode shapes for a given system is included in its right and left eigenvectors 
\cite{book:eigenvalue}. 
Given the eigenvectors of a system, an efficient measure of 
the shape of coupling between states and variables 
is provided through modal participation analysis \cite{arriaga:82_1}.
Consider system \eqref{eq:lin:ode}:
If 
$s_i$ is an eigenvalue of $\AS$
and $\reigv_i$, $\leigv_i$ are the associated eigenvectors, then the 
corresponding modal \acf{pf} is defined as the dimensionless number:\footnote{Definition \eqref{pf:eq:definition} assumes that the algebraic multiplicities of all eigenvalues equal the geometric ones. The reader interested in modal participation analysis of systems that do not satisfy this assumption is referred to \cite{dassios2020participation}.}
\begin{equation}\label{pf:eq:definition}
\PF
= \leigvel_{i,k}\, \reigvel_{k,i} \, ,
\end{equation}
where $\reigvel_{k,i}$ is the $k$-th row element of $\reigv_i$ and $\leigvel_{i,k}$ is the $k$-th column element of $\leigv_i$. The \ac{pf} in \eqref{pf:eq:definition} represents the relative contribution of the $i$-th mode $s_i$ in the response of the $k$-th state variable $x_k$.  Note that \acp{pf} can be collected to form the system's participation matrix $\PFmat$, as follows:  
\begin{equation}\label{eq:PFmat}
\PFmat = \leigvmat\T \circ \reigvmat  \, , 
\end{equation}
where $\circ$ denotes component-wise matrix multiplication; $\reigvmat$ is the modal matrix with the right eigenvectors as columns, and $\leigvmat$ is the modal matrix with the left eigenvectors as rows,~i.e.
$\reigvmat = 
 \begin{bmatrix}
     \reigv_1 & \reigv_2 & \ldots & \reigv_\nx  \\
 \end{bmatrix} 
 , \ 
  \leigvmat = 
 \begin{bmatrix}
     \leigv_1\T & \leigv_2\T & \ldots & \leigv_\nx\T  \\
 \end{bmatrix}\T$.

%

Now, consider a \ac{tdi} method and the associated approximated system \eqref{eq:difference}.
The modal participation matrix associated to \eqref{eq:difference} is then defined as follows \cite{book:eigenvalue}:
\begin{equation}\label{eq:PFmattdi}
\PFmattdi = \leigvmattdi\T \circ \reigvmattdi  \, ,
\end{equation}
with~
$ \reigvmattdi = 
\hspace{-0.5mm}
 \begin{bmatrix}
     \reigvtdi_1 & \reigvtdi_2 & \ldots & \reigvtdi_\ntdi  \\
 \end{bmatrix}$,  
  $\leigvmattdi = 
 \begin{bmatrix}
     \leigvtdi_1\T & \leigvtdi_2\T & \ldots & \leigvtdi_\ntdi\T  \\
 \end{bmatrix}\T $.
%
%
Note that matrix $\PFmattdi$ basically represents an approximation of the participation matrix $\PFmat$.  If $\PF$ is an element of $\PFmat$ and $\uppi$ is the corresponding element of $\PFmattdi$, then the quantity:
\begin{equation}\label{eq:p:relerr}
\errp = 100 {(|\uppi|-|\PF|)}/{|\PF|}  \, ,
\end{equation}
provides an estimate of the associated relative mode-shape deformation 
introduced by the \ac{tdi} method.

We note that metrics \eqref{eq:s:relerr} and \eqref{eq:p:relerr} are based on \ac{sssa} and thus they are technically valid around stationary solutions. Yet, the structure and stiffness of \eqref{eq:dae} as well as the properties of \ac{tdi} methods are features that tend to be “robust” and hence results provide also a tentative yet accurate estimate of deformation also for varying operating conditions.  For similar considerations we refer to the literature, e.g.,~\cite{arriaga:82_2,book:chow:13,tzounas2022tdistab}. 

\subsection{Deformation by Common Methods}
\label{sec:commuting}

In this section, we discuss the mode-shape deformation introduced by well-known \ac{tdi} methods used for the simulation of power system dynamics.
We first show that certain methods do not deform at all the mode-shapes of dynamic modes that are represented by \textit{non-degenerate} eigenvalues, i.e.~eigenvalues with algebraic multiplicity equal to 1.  To this aim, we provide the following result from linear algebra.

Consider two commuting matrices $\AS$ and $\Atdi$:
\begin{equation}\label{eq:commute}
    \AS \Atdi = \Atdi \AS \, .
\end{equation}
If $\reigv_i$ is a right eigenvector of $\AS$ corresponding to the non-degenerate eigenvalue $s_i$, then it is also an eigenvector of $\Atdi$.

\textit{Proof.}
The eigenvalue problem associated to $\AS$ is \eqref{eq:ep}, whereby substituting $\reigv_i,s_i$ and pre-multiplying by $\Atdi$ we get:
\begin{align}\nonumber
(s_i\Atdi -\Atdi\AS) \reigv_i &= 
(s_i\Atdi -\AS\Atdi) \reigv_i \\
&= (s_i\bfb I_\nx -\AS) \Atdi\reigv_i
=
\bfg 0_{\nx,1}  \, .
\label{eq:ep:comm}
\end{align}
Thus, $\Atdi \reigv_i$ is also a right eigenvector of $\AS$ associated to 
$s_i$ or, equivalently, $\Atdi\reigv_i$ is proportional to $\reigv_i$:
\begin{align}\label{eq:ep:comm:2}
(\lambda \bfb I_\nx -\Atdi) \reigv_i &= 
\bfg 0_{\nx,1}  \, ,
\end{align}
i.e.,~$\reigv_i$ is an eigenvector of $\Atdi$ associated to the eigenvalue $\lambda$.  Note that if $\Atdi$ represents a \ac{tdi} method, as is the case in this work, then $\lambda \equiv \zt_i$.  The reciprocal case of left eigenvectors can be derived similarly and thus, for a non-degenerate eigenvalue and a method whose matrix $\Atdi$ commutes with $\AS$, we have that $|\PF|=|\uppi|$ in \eqref{eq:p:relerr}, or, $\errp=0$. 

\subsubsection{Theta method}

Consider the Theta method \eqref{eq:theta}, for which $\Atdi$ is given by \eqref{eq:theta:pencil}.  To prove commutativity of $\Atdi$ and $\AS$, we start by considering the identity
$\AS - \ci\AS^2= \AS - \ci\AS^2$,
which can be equivalently 
%
%
rewritten as:
\begin{align}\label{eq:theta:rel3}
(\bfb I_\nx - \ci\AS) \AS 
  &= \AS (\bfb I_\nx - \ci\AS
) \, .
\end{align}
Left and right multiplication by 
$(\bfb I_\nx - \ci\AS)^{-1}$ yields:
\begin{align}\label{eq:theta:rel4}
\AS (\bfb I_\nx - \ci\AS)^{-1}  &= (\bfb I_\nx - \ci\AS)^{-1} \AS  \, .
\end{align}
Right multiplication of both sides of \eqref{eq:theta:rel4} by $\cii\AS$, $\cii \in \mathbb{R}$ gives:
\begin{align}\label{eq:theta:rel5}
\AS (\bfb I_\nx - \ci\AS
)^{-1}\cii\AS  &= (\bfb I_\nx - \ci\AS)^{-1}\cii\AS^2 \, .
\end{align}
Summing \eqref{eq:theta:rel4} and \eqref{eq:theta:rel5} and using $\ci=h(1-\theta)$, $\cii=h\theta$, leads to:
\begin{align}\nonumber
\AS [\bfb I_\nx -&h(1-\theta)\AS
]^{-1}(\bfb I_{\nx} + h\theta\AS) = \\
&=
\left[\bfb I_\nx - h(1-\theta)\AS
\right ]^{-1} (\bfb I_{\nx} + h\theta\AS)\AS \, ,
\label{eq:theta:rel6}
\end{align}
or, equivalently, to \eqref{eq:commute}.  The proof is complete.

\subsubsection{BEM and TM} 

They are special cases of the Theta method. 
For $\theta=0$, \ac{bem} is obtained and \eqref{eq:theta:rel6} becomes:
\begin{align}
\AS (\bfb I_\nx -h\AS
)^{-1} 
&=
(\bfb I_\nx - h\AS)^{-1} \AS \, .
\label{eq:bem:comm}
\end{align}
\ac{tm} is obtained for
$\theta=0.5$, in which case \eqref{eq:theta:rel6} becomes:
\begin{align}\nonumber
\AS (\bfb I_\nx -&0.5h\AS
)^{-1}(\bfb I_{\nx} + 0.5\theta\AS) = \\
&=
(\bfb I_\nx - 0.5h\AS)^{-1} (\bfb I_{\nx} + 0.5\theta\AS)\AS \, .
\label{eq:tm:comm}
\end{align}

\subsubsection{2S-DIRK}

We consider the \ac{dirk} proposed in \cite{2009noda} for the simulation of electromagnetic transients.  The method's first stage computes the solution at an intermediate point:
\begin{equation}
  \begin{aligned}
    \label{eq:2sdirk:1}
    \bfg \chi_{n+1} &= \bfg x_{n} + \dirka h \,
    \bfg f(\bfg \chi_{n+1}, \bfg \psi_{n+1}) \, , \\
     \bfg 0_{\ny,1} &= h \bfg g(\bfg \chi_{n+1},\bfg \psi_{n+1}) \, , \\
  \end{aligned}
\end{equation}
where $\dirka=1-1/\sqrt{2}$. Then, $\bfg \chi_{n+1}$ is used to calculate:
\begin{equation}
  \begin{aligned}
    \label{eq:2sdirk:2}
    \bfg \chi_{n} &= \dirkb \bfg x_{n} + (1-\dirkb) \, \bfg \chi_{n+1} \, , \quad \dirkb=-\sqrt{2} \, . \\
  \end{aligned}
\end{equation}
%
The final solution is obtained from the following equations:
\begin{equation}
  \begin{aligned}
    \label{eq:2sdirk:3}
    \bfg x_{n+1} &= \bfg \chi_{n}
    +\dirka h \, \bfg f(\bfg x_{n+1}, \bfg y_{n+1})
    \, ,\\
     \bfg 0_{\ny,1} &= h \bfg g(\bfg x_{n+1},\bfg y_{n+1}) \, .
  \end{aligned}
\end{equation}
Differentiating \eqref{eq:2sdirk:1}-\eqref{eq:2sdirk:3} at $(\bfg x_o,\bfg y_o)$ allows expressing the method in the form of \eqref{eq:difference}, where \cite{tzounas2022tdistab}:
\begin{equation}
  \begin{aligned}
    \label{eq:dirk:pencil}
    \Atdi = (\bfb I_\nx - \dirka h \AS )^{-1} ( \bfb I_\nx-\dirka \dirkb h \AS )
    (\bfb I_\nx-\dirka h \AS)^{-1}  .
  \end{aligned} 
\end{equation}
Matrices $\AS$, $\Atdi$ in \eqref{eq:dirk:pencil} are commuting. We omit the proof due to space constraints, but it can be easily constructed similarly to the Theta method starting from \eqref{eq:theta:rel4}, where in this case $\ci=\dirka h$.

Methods \textit{1)}-\textit{3)} above are implicit methods commonly used in a simultaneous-solution approach setup. 
We have shown that these methods do not deform the mode shape of dynamics represented by non-degenerate eigenvalues.  This is an important result since critical modes that dominate the dynamic response of power system models are typically represented by 
non-degenerate eigenvalues.  
The deformation introduced by these \ac{tdi} methods is further discussed through simulations in the case study presented in Section~\ref{sec:case}.

\subsubsection{Heun's Method (HM)}

We consider an 
element of the family of explicit Adams-Bashforth methods, namely \ac{hm}.  Variants of 
\ac{hm}
are commonly employed by software tools that adopt the partitioned-solution approach. 
In \ac{hm}, a predictor provides an initial estimate ($\bfg\xpc^{(0)}_{n+1}$) of $\bfg x_{n+1}$, as follows:
\begin{align}\label{eq:pred}
\bfg \xpc^{(0)}_{n+1} &= \bfg x_{n} + h 
\bfg f(\bfg x_{n},\bfg y_{n}) \, .
\end{align}

Then, accuracy of the current estimation is refined through corrector steps.  The $i$-th corrector step has the form:
\begin{equation}\label{eq:corr}
\bfg \xpc^{(i)}_{n+1} = \bfg x_{n}
+ 0.5h
\bfg f(\bfg x_{n},\bfg y_{n})
+  0.5h \bfg f (\bfg \xpc^{(i-1)}_{n+1}, 
\bfg y_{\rm n}) \, , 
\end{equation}
with $i \in \mathbb{N}^* : i\leq r$, where typically $r=1$ or $2$.  \ac{hm} needs to be combined with a way to 
deal with interfacing of algebraic variables \cite{machowski2020power}.  In \eqref{eq:corr},
such interfacing 
is achieved by extrapolation,
i.e.~$\bfg y_n$ is used instead of $\bfg y_{n+1}$ in the last term of the right-hand side of \eqref{eq:corr} 
\cite{tzounas:psastab,machowski2020power}. 
Then, $(\bfg x_{n+1},\bfg y_{n+1})$ is obtained from: 
%
\begin{align}
\label{eq:corr:xt}
\bfg x_{n+1} &= \bfg \xpc^{(r)}_{n+1} \, , \\
\bfg 0 &= h\bfg g(\bfg x_{n+1},\bfg y_{n+1}) \, .
\label{eq:pc:g}
\end{align}
Differentiation of \eqref{eq:corr}-\eqref{eq:pc:g} at $(\bfg x_o,\bfg y_o)$ allows expressing the method in the form of \eqref{eq:difference}, where:
\begin{equation}
  \begin{aligned}
    \label{eq:hm:pencil}
    \Atdi = \bfb I_\nx + h \sum_{j=0}^r \left(\frac{h}{2} \jac{f}{x}\right )^j  
    \hspace{-1mm} \AS  \, ,
  \end{aligned} 
\end{equation}
with $r \in \mathbb{N}^*$.  The proof of \eqref{eq:hm:pencil} can be found in \cite{tzounas:psastab}.
For $r\geq 1$, $\AS$ and $\Atdi$ in 
\eqref{eq:hm:pencil} do not commute, which implies that \ac{hm} is expected to deform the mode shapes of both degenerate and non-degenerate eigenvalues.  
If $r=0$, \ac{hm} reduces to the \ac{fem} and \eqref{eq:hm:pencil} yields $\Atdi = \bfb I_\nx + h \AS $.  In this case $\Atdi \AS = \AS \Atdi = \AS+h\AS^2$.  
Yet, \ac{fem} is known to show a poor performance, which from the viewpoint of this paper implies that the method gives rise to very large $\errs$ errors.  Thus, in the remainder of this work, \ac{fem} is not considered. 

\section{Case Study}
\label{sec:case}
In this section, we illustrate the proposed approach through simulations carried out based on the IEEE 39-bus test system.
The IEEE 39-bus system includes 10 \acp{sg} represented by a $4$-th order model, 34 lines, 12 transformers, and 19 loads.  All \acp{sg} are equipped with primary frequency and voltage
regulators, and power system stabilizers.  The system's static and dynamic data can be found in \cite{web:39bus}.
Simulations in this section are carried out using Dome~\cite{vancouver}.

\subsection{Eigenvalue Deformation}

\begin{figure*}[ht!]
     \centering
     \begin{subfigure}[b]{0.3\linewidth}
         \centering
         \includegraphics[width=\linewidth]{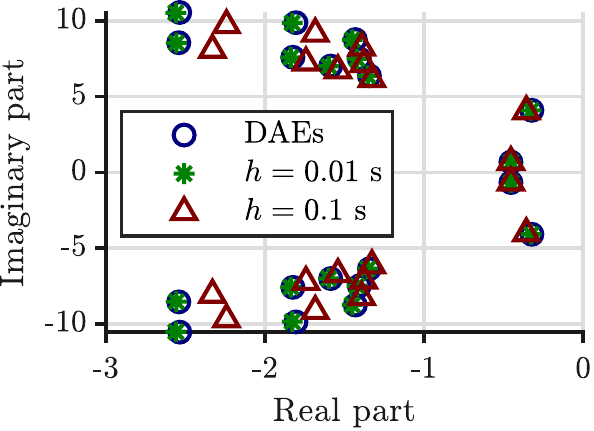}
         \caption{Theta ($\theta = 0.47$)}
         \label{fig:Eigs def ThetaM}
     \end{subfigure}
     \hfill
     \begin{subfigure}[b]{0.3\linewidth}
         \centering
         \includegraphics[width=\linewidth]{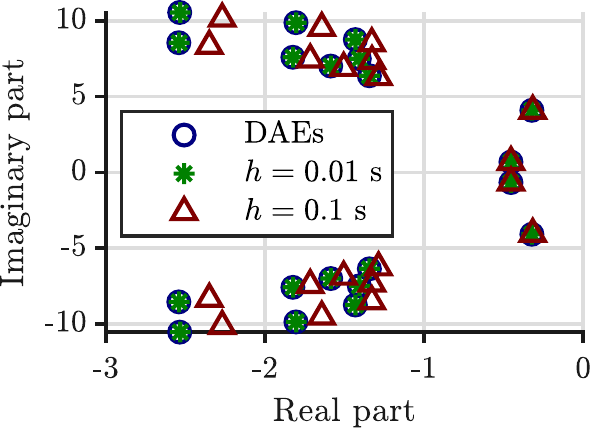}
         \caption{2S-DIRK}
         \label{fig:Eigs def 2S-DIRK}
     \end{subfigure}
     \hfill
     \begin{subfigure}[b]{0.3\linewidth}
         \centering
         \includegraphics[width=\linewidth]{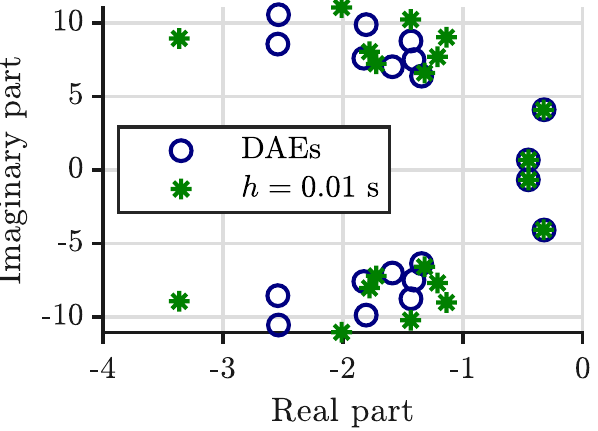}
         \caption{\ac{hm} ($r = 2$)}
         \label{fig:Eigs def HM}
     \end{subfigure}
     \vspace{-1mm}
        \caption{Eigenvalue deformation for Theta, \ac{dirk} and \ac{hm}.}
        \label{fig:Eigs def}
\end{figure*}

The eigenvalues of the DAE system obtained from \eqref{eq:ep} are compared to the ones of the associated problem \eqref{eq:ep:tdi} for Theta, \ac{dirk}, and \ac{hm}; numerical deformation is calculated as in \eqref{eq:s:relerr} for different time step sizes $h$.
Figure~\ref{fig:Eigs def} shows the spurious shift that 
these methods introduce to the rightmost eigenvalues of the system.  While \ac{hm} already presents considerable deviations from the exact system dynamics for $h = 0.01$~s,
Theta and \ac{dirk} have a good performance and notably deteriorate only for $h$ in the order of $10^{-1}$~s or higher. 

\begin{figure*}
\vspace{-1mm}
     \centering
     \begin{subfigure}[b]{0.3\linewidth}
         \centering
         \includegraphics[width=\linewidth]{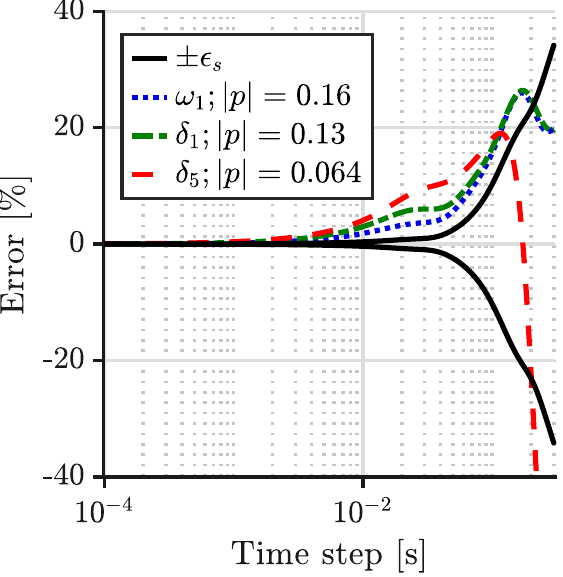}
         \caption{$s = -0.32 \pm 4.09\jj;$ $\zeta = 7.85$\%}
         \label{fig:Mode1}
     \end{subfigure}
     \hfill
     \begin{subfigure}[b]{0.3\linewidth}
         \centering
         \includegraphics[width=\linewidth]{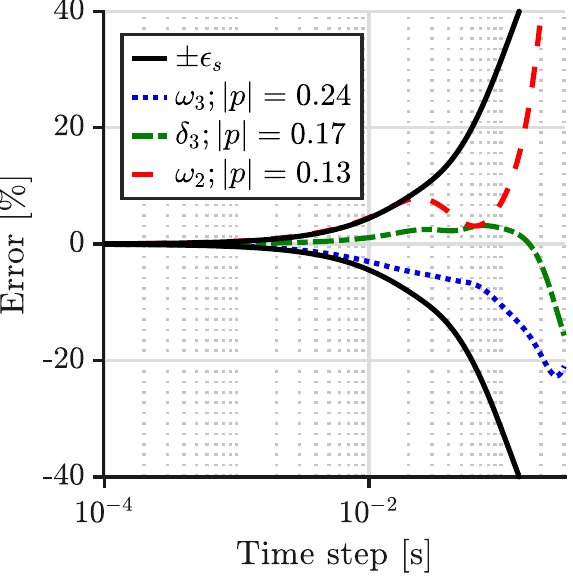}
         \caption{$s = -1.43 \pm 8.76\jj;$ $\zeta = 16.13$\%}
         \label{fig:Mode3}
     \end{subfigure}
     \hfill
     \begin{subfigure}[b]{0.3\linewidth}
         \centering
         \includegraphics[width=\linewidth]{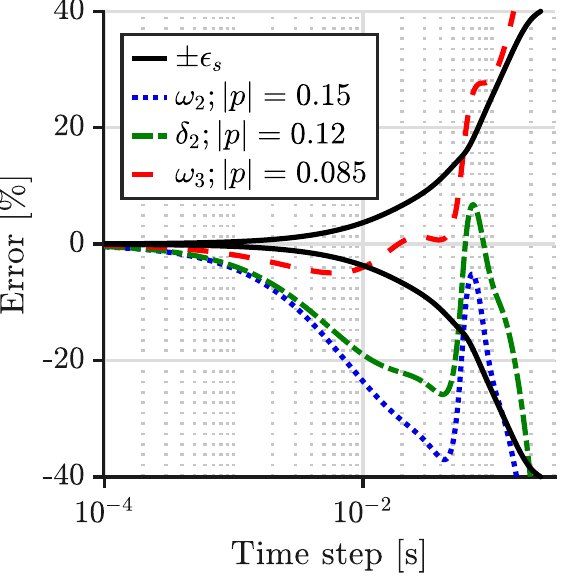}
         \caption{$s=-1.41 \pm 7.50\jj;$ $\zeta = 18.45$\%}
         \label{fig:Mode7}
     \end{subfigure}
     \vspace{-1mm}
        \caption{\ac{pf} deformation $\errp$ as a function of $h$ for 3 poorly damped modes 
        (with $\zeta$ denoting the damping ratio). 
        In each plot, the states with the highest \acp{pf} are represented, along with the associated eigenvalue deformation $\errs$.}
        \label{fig:PF def}   
\end{figure*}

\subsection{Deformation of Mode Shapes}


We focus on the deformation that \ac{tdi} methods introduce to the coupling between system states and variables.  To this aim,
the participation matrices 
$\PFmat$ and $\PFmattdi$
are calculated from \eqref{eq:PFmat} and \eqref{eq:PFmattdi}, respectively. Since different columns of $\PFmat$, $\PFmattdi$
refer to different modes, the columns of $\PFmattdi$ are sorted to pair correctly 
to the modes of the \ac{dae} system.  Moreover, the columns of both matrices are normalized so that for every eigenvalue the sum of all \acp{pf} 
is equal to 1.
%
%
Then, mode-shape deformation for each \ac{tdi} method is estimated from 
\eqref{eq:p:relerr}.  For Theta and \ac{dirk}, we find that $\errp=0$ for all non-degenerate eigenvalues, which is consistent with the discussion of
Section~\ref{sec:commuting}.  Furthermore,
for degenerate eigenvalues, large values of $\errp$ are observed in some cases.  Yet,
these cases are always associated with very low 
\acp{pf} ($|p| < 10^{-3}$).  Given that 
the behavior of a variable is 
largely defined by a small number
of highly participating modes (often by 1 or 2), the numerical impact of these cases on the system is negligible.

For \ac{hm}, 
significant values of 
$\errp$ are observed for 
both non-degenerate and degenerate 
eigenvalues. 
Figure~\ref{fig:PF def} shows,
for the most poorly damped electromechanical modes, how $\errp$ varies as a function of $h$ when 2 corrector steps are used. For the sake of comparison, 
the corresponding 
$\errs$ is included in each plot.
In Fig.~\ref{fig:PF def}, $\delta_i$, $\omega_i$ 
denote the rotor angle and speed, 
respectively, of the $i$-th \ac{sg}.
The deformation presents an irregular behavior but generally increases with the step size. 
Very small steps lead to good accuracy but also 
lead to a high computational burden.
Note also that for several modes and step sizes the maximum. 
$\errp$ is higher than 
$\errs$. 
For the most critical mode, for example
(Fig.~\ref{fig:Mode1}), $h=0.07$~s 
leads to $\errs<5$\% but also to a maximum 
$\errp$ of $\errp^{\max}>9$\%. 
Assuming for this mode a prescribed accuracy degree of $\errs,\errp^{\max}<5$\%, the maximum admissible time step is estimated at
$h^{\max}=0.012<0.07$~s. 
Another example is shown in Fig.~\ref{fig:Mode7}, where
$h^{\max}=1$~ms is needed to maintain 
$\errp^{\max}<5$\%, 
although $\errs$
is low even for $h=0.01$~s.
The above discussion highlights the relevance of evaluating both metrics in a numerical analysis.

\vspace{-2mm}

Table~\ref{tab:maxh} reports, for different methods, the maximum admissible time step 
$h^{\max}$ in 4 scenarios:
(i) $\errs<5$\%,
(ii) $\errp<5$\%, 
(iii) $\errp<10$\%,
and
(iv) $\errs,\errp<5$\%. 
In all cases, $h^{\max}$ is obtained considering the 5 most critical eigenvalues -- which are all non-degenerate -- and for each eigenvalue, the 3 largest \acp{pf}.
Results show that
$h^\mathrm{max}$ for implicit
methods is about an order of magnitude larger than \ac{hm} and is not impacted by mode-shape deformation.  On the other hand, 
$h^\mathrm{max}$ for \ac{hm} is largely impacted by the selected $\errp$ threshold. 
We note that in practice $\errp^{\max}$, $\errs^{\max}$ can be setup for any prescribed requirements.  In this regard, a relevant question that is worth further study is how to best tune $\errp^{\max}$, $\errs^{\max}$.
A good starting point in this direction can be the analytical solution of the linearized system, which depends linearly on eigenvectors (and thus on mode shapes), but exponentially on system eigenvalues, which directs that  
$\errp^{\max}>\errs^{\max}$.

\begin{table}[!t]
    \centering
    \footnotesize 
    \renewcommand{\arraystretch}{1.15}
    \caption{$h^{\max}$ 
    estimated for the 5 most 
    critical modes and 3 highest \acp{pf}. 
 Theta and \ac{dirk}
 have $\errp=0$.}
    \label{tab:maxh}
    \begin{tabular}{l|ccccc}
    \toprule
      \multirow{2}{*}{Accuracy} 
      & \multicolumn{4}{c}{Method}
      \\
      &  
      {Theta}
      & 
      {\ac{dirk}}
       & 
       \ac{hm} ($r=1$) & \ac{hm} ($r=2$)\\ 
      \midrule
      $\errs<5$\%
    & 
    0.080 & 0.115 & 
    0.0087 & 0.0098 
    \\
    $\errp<5$\%
    & $\infty$ & $\infty$ & 
    0.0012 & 0.0012  \\
    $\errp<10$\%
    & $\infty$ & $\infty$ & 
    0.0026 & 0.0027  \\
    $\errs,\errp<5$\%
    & 0.080 & 0.115 & 
    0.0012 & 0.0012  
    \\ \bottomrule
    \end{tabular}
\end{table}

\subsection{Modified System with DERs and AGC}

\begin{figure}
    \centering
    \includegraphics[width=\linewidth]{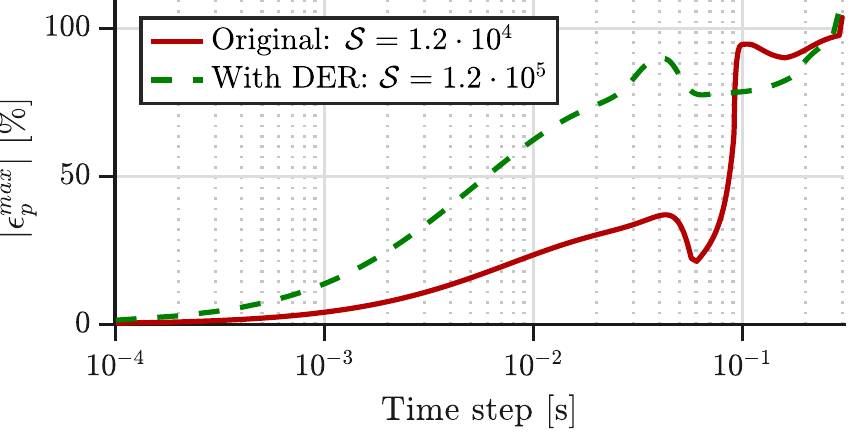}
    \caption{Effect of stiffness on \ac{hm} mode-shape deformation. The curves show $\errp^{\max}$ for the original and modified systems.}
    \label{fig:stiff}
\end{figure}

In this section, the test system is modified as follows.    
\acp{sg} at buses 32, 33, 34 and 35 are replaced by aggregated, converter-based \ac{der} models.  Each \ac{der} synchronizes to the grid through a synchronous reference frame \ac{pll} and provides primary frequency and voltage support by regulating, at the point of connection, the $d$ and $q$ axis current components, respectively, in the $dq$ reference frame.  Moreover, \acp{sg} 
are assumed to provide secondary frequency support 
through an \ac{agc} scheme modeled as an integral regulator.  The \ac{der} \ac{pll} dynamics are faster than the fastest dynamics 
of the original system, while the \ac{agc} dynamics are slower than the slowest dynamics of the original system. As a consequence, 
 the system's stiffness ratio, defined as in \eqref{eq:stiffness}, increases by an order of magnitude (from $1.2 \cdot 10^4$ to $1.2 \cdot 10^5$).
Figure~\ref{fig:stiff} shows, as a function 
of $h$, the relative error
$\errp^{\max}$ for the 3 most participating states of each of the 5 least damped modes.  It is seen that increasing the system's stiffness results in a higher distortion, particularly for small time steps.  

    
\section{Conclusion}
\label{sec:conclusion}

This paper studies the numerical deformation that \ac{tdi} methods cause to the mode shapes of power system \ac{dae} models.  It is shown that, owing to matrix commutativity properties, 
common implicit methods, such as 
Theta and \ac{dirk}, do not deform 
the mode shape of dynamics represented by non-degenerate eigenvalues.  
Moreover, the well-known \ac{hm} is employed to illustrate through simulations the effect on mode shapes for the case that commutativity properties do not hold.  In future work, we will 
employ the proposed approach to 
study the numerical robustness
of state-of-art converter-based models 
integrated with explicit \ac{tdi} methods.

\bibliographystyle{IEEEtran}
\bibliography{references}

\begin{thebibliography}{10}
\providecommand{\url}[1]{#1}
\csname url@samestyle\endcsname
\providecommand{\newblock}{\relax}
\providecommand{\bibinfo}[2]{#2}
\providecommand{\BIBentrySTDinterwordspacing}{\spaceskip=0pt\relax}
\providecommand{\BIBentryALTinterwordstretchfactor}{4}
\providecommand{\BIBentryALTinterwordspacing}{\spaceskip=\fontdimen2\font plus
\BIBentryALTinterwordstretchfactor\fontdimen3\font minus
  \fontdimen4\font\relax}
\providecommand{\BIBforeignlanguage}[2]{{%
\expandafter\ifx\csname l@#1\endcsname\relax
\typeout{** WARNING: IEEEtran.bst: No hyphenation pattern has been}%
\typeout{** loaded for the language `#1'. Using the pattern for}%
\typeout{** the default language instead.}%
\else
\language=\csname l@#1\endcsname
\fi
#2}}
\providecommand{\BIBdecl}{\relax}
\BIBdecl

\bibitem{kundur:94}
P.~Kundur, \emph{{Power System Stability and Control}}.\hskip 1em plus 0.5em
  minus 0.4em\relax New York: Mc-Grall Hill, 1994.

\bibitem{machowski2020power}
J.~Machowski, Z.~Lubosny, J.~W. Bialek, and J.~R. Bumby, \emph{Power system
  dynamics: stability and control}.\hskip 1em plus 0.5em minus 0.4em\relax John
  Wiley \& Sons, 2020.

\bibitem{1995paserba}
J.~Sanchez-Gasca, R.~D'Aquila, W.~Price, and J.~Paserba, ``Variable time step,
  implicit integration for extended-term power system dynamic simulation,'' in
  \emph{Proceedings of Power Industry Computer Applications Conference}, 1995,
  pp. 183--189.

\bibitem{powerfactory}
{DIgSILENT Power System Solutions}, ``{DIgSILENT PowerFactory},''
  \href{https://www.digsilent.de/powerfactory.html}{digsilent.de/powerfactory}.

\bibitem{stott:1979}
B.~{Stott}, ``Power system dynamic response calculations,'' \emph{Proceedings
  of the {IEEE}}, vol.~67, no.~2, pp. 219--241, Feb. 1979.

\bibitem{psse_pag2}
{PSS/E 33.0}, \emph{Program Application Guide Volume 2}.\hskip 1em plus 0.5em
  minus 0.4em\relax Siemens, 2011.

\bibitem{ge:pslf}
{General Electric Energy Consulting}, ``{General Electric (GE) PSLF},''
  \href{https://www.geenergyconsulting.com/practice-area/software-products/pslf}{geenergyconsulting.com/practice-area/software-products/pslf}.

\bibitem{ramasubramanian2020positive}
D.~Ramasubramanian, W.~Wang, P.~Pourbeik, E.~Farantatos, A.~Gaikwad, S.~Soni,
  and V.~Chadliev, ``Positive sequence voltage source converter mathematical
  model for use in low short circuit systems,'' \emph{IET Generation,
  Transmission \& Distribution}, vol.~14, no.~1, pp. 87--97, 2020.

\bibitem{ramasubramanian2022parameterization}
D.~Ramasubramanian, X.~Wang, S.~Goyal, M.~Dewadasa, Y.~Li, R.~O’Keefe, and
  P.~Mayer, ``Parameterization of generic positive sequence models to represent
  behavior of inverter based resources in low short circuit scenarios,''
  \emph{Electric Power Systems Research}, vol. 213, p. 108616, 2022.

\bibitem{tzounas2022tdistab}
G.~Tzounas, I.~Dassios, and F.~Milano, ``Small-signal stability analysis of
  numerical integration methods,'' \emph{IEEE Transactions on Power Systems},
  vol.~37, no.~6, pp. 4796--4806, Nov. 2022.

\bibitem{tzounas2022tdistabd}
------, ``Small-signal stability analysis of implicit integration methods for
  power systems with time delays,'' \emph{Electric Power System Research}, vol.
  211, no. 108266, Oct. 2022.

\bibitem{tzounas:psastab}
G.~Tzounas and G.~Hug, ``Unified numerical stability and accuracy analysis of
  the partitioned-solution approach,'' \emph{submitted to IEEE Transactions on
  Power Systems}, 2022, under review. Available at:
  {\href{https://n.ethz.ch/~gtzounas/pap/psastab.pdf}{n.ethz.ch/$\sim$gtzounas/pap/psastab.pdf}}.

\bibitem{2022ampsas}
F.~Milano, M.~Liu, M.~A.~A. Murad, G.~M. J{\'o}nsd{\'o}ttir, G.~Tzounas,
  M.~Adeen, {\'A}.~Ortega, and I.~Dassios, ``Power system modelling as
  stochastic functional hybrid differential-algebraic equations,'' \emph{IET
  Smart Grid}, vol.~5, no.~5, pp. 309--331, Oct. 2022.

\bibitem{book:eigenvalue}
F.~Milano, I.~Dassios, M.~Liu, and G.~Tzounas, \emph{Eigenvalue Problems in
  Power Systems}.\hskip 1em plus 0.5em minus 0.4em\relax CRC Press, Taylor \&
  Francis Group, 2020.

\bibitem{arriaga:82_1}
I.~J. P{\'e}rez-Arriaga, G.~C. Verghese, and F.~C. Schweppe, ``Selective modal
  analysis with applications to electric power systems, part i: Heuristic
  introduction,'' \emph{IEEE Transactions on Power Apparatus and Systems}, vol.
  PAS-101, no.~9, pp. 3117--3125, Sep. 1982.

\bibitem{dassios2020participation}
I.~Dassios, G.~Tzounas, and F.~Milano, ``Participation factors for singular
  systems of differential equations,'' \emph{Circuits, Systems, and Signal
  Processing}, vol.~39, no.~1, pp. 83--110, 2020.

\bibitem{arriaga:82_2}
G.~C. Verghese, I.~J. P{\'e}rez-Arriaga, and F.~C. Schweppe, ``Selective modal
  analysis with applications to electric power systems, part ii: the dynamic
  stability problem,'' \emph{IEEE Transactions on Power Apparatus and Systems},
  vol. PAS-101, no.~9, pp. 3126--3134, Sep. 1982.

\bibitem{book:chow:13}
J.~H. Chow, \emph{{Power System Coherency and Model Reduction}}, ser. Power
  Electronics and Power Systems 94.\hskip 1em plus 0.5em minus 0.4em\relax New
  York: Springer-Verlag, 2013.

\bibitem{2009noda}
T.~Noda, K.~Takenaka, and T.~Inoue, ``Numerical integration by the 2-stage
  diagonally implicit {Runge-Kutta} method for electromagnetic transient
  simulations,'' \emph{IEEE Transactions on Power Delivery}, vol.~24, no.~1,
  pp. 390--399, 2009.

\bibitem{web:39bus}
{Illinois Center for a Smarter Electric Grid (ICSEG)}, ``{IEEE 39-Bus
  System},''
  \href{http://publish.illinois.edu/smartergrid/ieee-39-bus-system/}{publish.illinois.edu/smartergrid/ieee-39-bus-system/}.

\bibitem{vancouver}
F.~Milano, ``{A Python-based software tool for power system analysis},'' in
  \emph{Proceedings of the IEEE PES General Meeting}, Jul. 2013.

\end{thebibliography}
  

\end{document}